\documentclass[a4]{article}

\usepackage{amsmath}
   \usepackage{paralist}
   
  \usepackage[colorlinks=true]{hyperref}
  \hypersetup{urlcolor=blue, citecolor=red}

\usepackage{amssymb} 
\usepackage{amsthm} 
\usepackage[english]{babel}
\usepackage{bbm}

\numberwithin{equation}{section}

\theoremstyle{definition}

\newtheorem{remark}{Remark}
\newtheorem*{notation}{Notation}

\newcommand{\R} {\ensuremath {\mathbb{R}}}

\newcommand{\calI} {\ensuremath {\mathcal{I}}}
\newcommand{\calJ} {\ensuremath {\mathcal{J}}}
\newcommand{\calF} {\ensuremath {\mathcal{F}}}
\newcommand{\calH} {\ensuremath {\mathcal{H}}}
\newcommand{\calY} {\ensuremath {\mathcal{Y}}}
\newcommand{\calK} {\ensuremath {\mathcal{K}}}
\newcommand{\calS} {\ensuremath {\mathcal{S}}}

\newcommand{\sinc} {\mathrm{sinc}}
\newcommand{\dd}  {\ensuremath{\,\mathrm{d}}}
\newcommand{\indicator}[1]{\mathbbm{1}_{\left\{ {#1} \right\} }}

\DeclareMathAlphabet{\itbf}{OML}{cmm}{b}{it}

\def\by{{\itbf y}}
\def\bx{{\itbf x}}

\def\bd{{\itbf d}}
\def\bu{{\itbf u}}
\def\bv{{\itbf v}}
\def\bw{{\itbf w}}

      \title{High frequency analysis of imaging with noise blending}

\author{ {\scshape Ennio Fedrizzi} \\ 
LPMA, Universit\'e Paris 7 Diderot, Sorbonne Paris Cit\'e\\
75013 Paris, France.\\
{\small fedrizzi@math.univ-paris-diderot.fr}
}

\date{}

\begin{document}

\maketitle




\begin{abstract}
We consider sensor array imaging for simultaneous noise blended sources. We study a migration imaging functional and we analyze its sensitivity to singular perturbations of the speed of propagation of the medium. We consider two kinds of random sources: randomly delayed pulses and stationary random processes, and three possible kinds of perturbations. Using high frequency analysis we prove the statistical stability (with respect to the realization of the noise blending) of the scheme and obtain quantitative results on the image contrast provided by the imaging functional, which strongly depends on the type of perturbations.
\end{abstract}  \bigskip 

\bigskip 
\bigskip

\section{Introduction}

An amazing fact in the analysis of imaging functionals, that has been recently pointed out and is currently under investigation, is that specific kinds of noise can improve the image quality or drastically reduce computational costs associated to the evaluation of the imaging functional.
A strong motivation is the observation that time reversal refocusing is enhanced when the medium is randomly scattering. A time reversal experiment is based on the use of a special device called a time reversal mirror (TRM), which is an array of transducers that can be used as receivers and transmitters. A typical time reversal experiment consists in two steps. In a first step, a point source emits a short pulse that propagates through the medium and is recorded by the TRM used as an array of receivers. In a second step, the recorded signals are time reversed and reemitted by the TRM used as an array of transmitters. The waves then refocus on the original point source location. The striking observation is that refocusing is enhanced when the medium is randomly heterogeneous and scattering compared to the situation in which the medium is homogeneous. Moreover, the refocused pulse is statistically stable in the sense that it does not depend on the realization of the random medium, but only on the statistical distributions of the fluctuations of the random medium \cite{FGPS}.

In the context of sensor array imaging, a similar technique is employed. The typical experiment still consists in two steps. The first step is the experimental data acquisition: a point source emits a wave into the medium, the wave is reflected by the singularities in the propagation speed of the medium and is recorded by an array of receivers. The second step is a numerical processing of the recorded data: the recorded signals are time reversed and resent in a numerical simulation into a model medium to locate the singularities in the propagation speed. However, in contrast with physical time reversal, the fictitious medium employed in the imaging process cannot in general capture complex heterogeneities of the original medium. Therefore, research has focused on other approaches to improve the imaging process,  such as the use of random sources \cite{DFGS12}, \cite{GP10}, \cite{HSH08}, \cite{SW11}, \cite{VB11}, \cite{WNT12}.

The classical approach to the imaging problem consists in performing a large number of experiments sounding each time a different source. For each experiment, the signal recorded at each receiver is stored and time reversed. This produced the data matrix. For each experiment, the time reversed data are reemitted (numerically) into the model medium in a new simulation, and the images obtained by each simulation are stacked. 
While this method provides a very good \textquotedblleft image\textquotedblright of the medium, it involves gathering, storing and processing huge amounts of data \cite{Be09}, \cite{MDB11}.

In the \emph{noise blending} approach noisy sources are used and they are all sounded simultaneously in one experiment. In this case, the time reversed recorded signals from the physical experiment are stored in a data vector and resent simultaneously into the model medium in a single simulation. This approach allows for considerable savings in both the data gathering, storing and processing stages. However, special care must be put in the choice of the noisy sources in order to ensure that the cross talk terms are very small and do not compromise image quality. This technique can be successfully applied also to the time reversal approach and its analysis bears similarities with techniques for passive imaging, which exploits ambient noise sources to recover travel times from correlations in between recordings at different stations \cite{GP09}. Applications of these different techniques are being investigated in different fields, from seismology \cite{SC06} \cite{LM06} \cite{SC05} \cite{GS08}, to volcano monitoring \cite{SR06} \cite{BS08} \cite{BS07}, to petroleum prospecting \cite{CG06} and medicine \cite{kidney}.

As remarked, the crucial step for the noise blending approach lies in the choice of the noisy sources. In \cite{DFGS12} it was suggested the use of stationary random Gaussian sources or of randomly delayed source pulses: for these choices of random sources, fourth moment computations show that the algorithm is statistically stable. In the present work we develop some further investigations on this setting and show how to obtain quantitative results on image quality and statistical stability of this algorithm in the high frequency regime, when the goal is to image singular perturbations in the speed of propagation. The result strongly depends on the type of perturbation. We therefore consider three types of perturbations, supported respectively on small balls, thin tubes and thin discs. With a slight abuse of notation we will refer to them as point, line and surface singularities, as they can be thought of as approximations of singular perturbations of the velocity of propagation supported on subspaces of lower dimension.

For each kind of perturbation we first analyze the average image contrast seen between the center of the perturbation and a point far from it: this will provide an hint on the level of difficulty to correctly image these perturbations. An even more interesting result follows: it concerns the quantitative analysis of the statistical stability of this functional, providing the typical contrast seen for the three perturbations. The question of stability of the imaging functional has already been addressed in \cite{DFGS12}: no quantitative analysis was carried out there, but it was shown that a condition for the statistical stability is that the recording time interval $T$ must be large. With a careful analysis of fluctuations produced by the random sources we show that, in the high frequency regime, the typical contrast is actually much better than just of order $\sqrt{T}$; but again the exact order of amplitude of the contrast depends on the shape of the perturbation. Point singularities are easy to observe, while surface type perturbations are the hardest to locate.

The paper is organized as follows. A short presentation of the model and the imaging functional used will complete this introductory section. In section \ref{sec:E} we analyze the average (with respect to the realization of the random time delays used in the blending process) sensitivity of the method to the three types of singularities. In section \ref{sec:oscillation} we study fluctuations due to the stochastic nature of the result, and from the analysis of the typical behavior we obtain conditions ensuring the possibility to accurately image the perturbations. Finally, all the results obtained are collected and discussed in section \ref{sec:end}.\\

\begin{notation} We use boldfaced characters to denote vectors: for example $x\in\R$, but $\bx = (x,y,z) \in \nolinebreak \R^3$. $B_r$ will denote the ball of radius $r$ in $\R^d$ for $d$ equals 2 or 3: $B_r=\{ \bx \in \R^d |\, |\bx| \le r\}$. The notation $A=O(\varepsilon)$ means that the quantity $A$ is \emph{exactly} of order $\varepsilon$; in order to say that it is of order $\varepsilon$ or smaller we will write $A\le O(\varepsilon)$. 
\end{notation}

\subsection{The wave equation}
We consider the solution $u$ of the wave equation in a three-dimensional inhomogeneous medium
\begin{equation}
\label{eq:wave intro}
\frac{1}{c^2(\bx)} \frac{\partial^2 u}{\partial t^2} - \Delta_\bx u = {n}(t,\bx)\ ,
\end{equation}
where $c(\bx)$ is the velocity of propagation of waves in the medium and $\bx=(x,y,z)\in\R^3$. We rewrite the velocity in the form
\begin{equation*}
c^{-2}(\bx) = c_0^{-2}(\bx) + \delta c^{-2} (\bx) \ ,
\end{equation*}
where $ {c_0(\bx)} $ is the known {\bf smooth} background velocity (for simplicity we assume it to be constant) and $\delta c^{-2}(\bx)$ is the velocity perturbation that we want to estimate, whose spatial support is contained in some domain $\Omega \subset \R^3$. We will detail these perturbations in the following. To simplify the geometry of the model, we shall take $\Omega=B_R$ to be a ball with a large radius $R$. We also assume to have $N_s$ point sources located at points $(\by_s)_{s =1,\ldots,N_s}$ laying on the surface $\partial B_R$. They can either emit (almost) simultaneously the same short pulse waveform, but randomly delayed in time, or independent stationary random signals. Following \cite{DFGS12}, we will refer to the first case as (\emph{noise}) \emph{blended} sources and to the second as \emph{stationary random} sources.

For the first case, the source term $n(t,\bx)$ is of the form
\begin{equation*}  
n(t,\bx) = \sum_{s=1}^{N_s} f(t-\tau_s) \delta(\bx-\by_s) \ .
\end{equation*}
The (short) pulse function $({f}(t) )_{t \in \R}$ is deterministic. Its carrier frequency is $\omega_0$ and its bandwidth is $b$.
The time delays $(\tau_s)_{s =1,\ldots,N_s}$ are zero-mean independent and identically distributed random variables with probability density function $p_\tau(t)$.   

For the second case, the source term is given by
\begin{equation*}
n(t,\bx) =  \sum_{s=1}^{N_s}  n_s(t)  \delta(\bx-\by_s) \ ,
\end{equation*}
where the random functions $({n}_s(t) )_{t \in \R}$, $s=1,\ldots,N_s$, are independent, zero-mean, stationary Gaussian processes with autocorrelation function
\begin{equation*} 
\left< {n}_s(t_1 ) {n}_{s'}(t_2 ) \right> = \delta_{ss'}  {F} (t_2-t_1) \, .
\end{equation*}\\

The direct and inverse problems can be formulated in terms of the background Green's function, the fundamental solution of the wave equation (\ref{eq:wave intro}). For an homogeneous medium with constant background velocity $c_0$, in the Fourier domain the Green's function is given by
\begin{equation*}
\hat G(\omega,\bx_1, \bx_2) = \frac{1}{4\pi |\bx_1-\bx_2|} \exp \Big( i \frac{\omega}{c_0} |\bx_1-\bx_2| \Big) \ .
\end{equation*}
Here, the Fourier transform of a function $f(t)$ is defined by
\begin{equation*}
\hat{f}(\omega ) =\int f(t) e^{i \omega t} \dd t \, .
\end{equation*}

\subsection{Direct and inverse problems}\label{sec:direct-inverse Pb}
We introduce here the scattering operator, that is, the mapping from velocity perturbations to the data, in the Born
approximation \cite{bleistein01}.

We assume that sources (located at points $\by_s$, $s=1, ... ,N_s$)  are disposed on the surface of the ball $B_R$ containing the perturbations, and are dense enough (ideally, closer than half of the central wavelength) so that a continuum approximation can be used. Signals are observed at the passive sensor array $(\bx_r)_{r=1,\ldots,N_r}$ for some large time interval $[-T/2,T/2]$. For noise blended sources, the recording time $T$ should be much larger than the typical travel time, so as to guarantee that the backscattered signals are completely recorded. For stationary random sources it must be taken much larger that the inverse of the bandwidth of the noise sources (i.e. much larger than the decoherence time).

The recorded data consists of the vector ${\it \bd}(t)=(d(t,\bx_r ))_{r =  1,\ldots,N_r}$ of the signals recorded by $\bx_r$ for $t\in [-T/2,T/2]$. 
These data are modeled by the scattering operator
${\mathcal F} :\ \big(\delta c^{-2}(\bx) \big)_{\bx \in B_R} \to \big({\it \bd}(t) \big)_{t\in [-T/2,T/2]}$, where
\begin{equation*} 
\big({\mathcal F} \delta c^{-2} \big)(t,\bx_r)  =\int_{\Omega} \hspace*{-0.03in} Q(t,\bx_r,\bx) \delta c^{-2}(\bx) \dd\bx \, ,
\end{equation*}
\begin{equation}\label{eq:def Q}   
\hspace{-0.1mm} 
Q(t,\bx_r,\bx) = - \frac{\partial^2}{\partial t^2}   \hspace*{-0.03in}   \iiint \hspace*{-0.03in} G (t_1 ,\bx_r,\bx) G(t_2,\bx,\by)  {n}(t-t_1-t_2,\by)  \dd t_1 \dd t_2 \dd \by \, . \hspace{-0.2mm}
\end{equation}
In the Fourier domain
\begin{equation}\label{eq:def Q-fourier}  
\hat{Q}( \omega ,\bx_r,\bx) = \omega^2  
   \int  \hat{G}(\omega,\bx_r,\bx) \hat{G}(\omega,\bx,\by)\hat{{n}}(\omega,\by) \dd\by \, .
\end{equation}
This is the formulation of the direct problem: the expression of the data set in terms of the velocity perturbation.\\

The imaging problem (inverse problem) aims at inverting the map $\calF$ in order to reconstruct the velocity perturbation $\big(\delta c^{-2}(\bx) \big)_{\bx \in B_R}$ from the data set $\big({\it \bd}(t) \big)_{t\in [-T/2,T/2]}$. The usual (least--square) approach would consist in applying the operator $(\calF^*\calF)^{-1} \calF^*$ to the data set ${\it \bd}$, where the adjoint of the scattering operator $\calF$ is given by
\begin{equation} \label{eq:def F*}
\big( {\mathcal F}^*  {\itbf d} \big)(\bx)
= \sum_{r=1}^{N_r} \int_{-\frac{T}{2}}^{\frac{T}{2}}  Q(t,\bx_r,\bx)  d(t,\bx_r) \dd t  \, .
\end{equation}
However, the full least--square inversion is in practice too complicated and the normal operator ${\mathcal F}^*{\mathcal F} $ is usually dropped in the inversion process. In \cite{DFGS12} it was shown that for $T$ large the normal operator is statistically stable (i.e. its fluctuations  are smaller than its expectation) and that its statistical average is close to the identity operator (more precisely, the kernel $\langle \calF^*\calF \rangle (\bx,\bx')$ concentrates near the diagonal $\bx=\bx'$), proving that this procedure can indeed provide a reasonable estimate of the velocity perturbation.

In the following sections we perform a detailed analysis of this imaging functional in the high frequency regime, obtaining quantitative results on the statistical stability of the method for different types of perturbations.

Since the main advantage of our approach lays in the drastic reduction in storage and computational costs, let us stress that the computation of the adjoint operator $\calF^*$ can be done quite easily.
\medskip
\begin{remark} The adjoint operator $\calF^* $ can be computed solving only two wave equations as follows. First, compute the wave $u(t,\bx)$ emitted by the original source, which is to say solve the wave equation with source $n(t,\bx)$ and background velocity $c_0$:
\begin{equation*}
u(t,\bx) = \iint G(t_1,\bx,\by) n(t-t_1,\by) \dd \by \dd t_1 \ .
\end{equation*}
 Second, compute the anti causal solution $v(t,\bx)$ of the wave equation with source term \linebreak $\sum_{r} \delta(\by-\bx_r) \partial_t^2 d(t,\bx_r)$:
\begin{align*}
v(t,\bx) &= \sum_{r=1}^{N_r} \int G(t_2-t,\bx,\bx_r) \, \partial^2_{t_2} d(t_2,\bx_r) \dd t_2 = \sum_{r=1}^{N_r} \int G(t_2,\bx,\bx_r) \, \partial^2_t d(t_2+t,\bx_r) \dd t_2  \, .
\end{align*}
Correlating the two wave solutions at a point $\bx$ in the search window produces the imaging functional
\begin{align*}
\calI(\bx)& = -\int u(t,\bx) v(t,\bx) \dd t \\
&=  \sum_{r=1}^{N_r}  - \partial^2_{t} \int \iiint G(t_2,\bx,\bx_r) G(t_1,\bx,\by) n(t-t_1-t_2,\by) d(t,\bx_r)  \dd t_1 \dd t_2 \dd \by \, \dd t
\end{align*}
and using the definition (\ref{eq:def Q}) we obtain 
\begin{equation*}
\calI(\bx) = \big( {\mathcal F}^*  {\itbf d} \big)(\bx) \ .
\end{equation*}
\end{remark}

\subsection{Analysis of the imaging functional}\label{sec:imaging func}
Let us start by the analysis of a kernel which will appear in the imaging functional we have to study. 

In our model, sources and receivers are located on the surface of the ball $B_R$ containing the perturbations and are dense enough so that a continuum approximation can be used. We can therefore rewrite the kernel
\begin{align*}
\calK_\omega(\bx,\by) &:= \sum_{r=1}^{N_r} \overline{\hat{G}}\big(\omega, \bx_r, \bx ) \hat{G}\big(\omega, \bx_r, \by  \big) \simeq \int_{\partial B_R} \rho(\bx_r) \overline{\hat{G}}\big(\omega, \bx_r, \bx ) \hat{G}\big(\omega, \bx_r, \by  \big) \dd \sigma(\bx_r) \, ,
\end{align*}
where $\rho(\bx_r)$ is the surface density of receivers or sources. We also assume that perturbations are located near the center of the ball; this means that the medium is homogeneous outside of a ball $B_r$ of radius $r\ll R$. 
Then, we have the approximate identity \cite[Proposition 4.3]{GP09}
\begin{equation*}
\frac{2 i \omega}{c_0} \int_{\partial B_R} \overline{\hat{G}}\big(\omega, \bx_r, \bx ) \hat{G}\big(\omega, \bx_r, \by  \big) \dd \sigma(\bx_r) 
\simeq \hat{G}\big(\omega, \bx, \by  \big) - \overline{\hat{G}}\big(\omega, \bx, \by ) \, ,
\end{equation*}
which follows from Green's identity and the Sommerfeld radiation condition. This result can be viewed as a version of the Helmholtz-Kirchhoff integral theorem. Using the function $sinc(x) = sin(x)/x$, the right hand side of the above equation can be rewritten as
\begin{equation*}
2 i \, \Im \bigg(  \hat{G}\big(\omega, \bx, \by  \big) \bigg) = \frac{2i}{4\pi} \, \frac{\omega}{c_0} \sinc\Big(\frac{\omega}{c_0} |\bx-\by|\Big) \,  
\end{equation*}
and assuming that receivers and sources have constant density on the surface $\partial B_R$ we obtain
\begin{equation*}
\calK_\omega(\bx,\by) = \frac{1}{4\pi} \sinc\Big(\frac{\omega}{c_0} |\bx-\by|\Big) \, .
\end{equation*}

We return to the analysis of the imaging functional, which estimates the velocity perturbations, given by (\ref{eq:def F*}). Using equation (\ref{eq:def Q-fourier}) and the definition of the vector ${\it \bd}(t)$ of recorded signals, we can write the imaging functional as
\begin{align} 
\widehat{\delta c} {}^{-2} (\bx) &= \big( \calF^* {\itbf d} \big) (\bx) = \frac{1}{2\pi} \sum_{r=1}^{N_r} 
      \iint  \omega^2 {\hat{G}}(\omega, \bx,\bx_r) {\hat{G}}( \omega, \bx,\by) 
  \hat{n}(\omega,\by)\overline{\hat{d}} (\omega, \bx_r ) \, \dd \by \dd \omega \,  \nonumber\\
  &=\int_{B_R} \calF^*\calF (\bx,\bx') \delta c^{-2} (\bx') \dd \bx'   = \int_{B_R} \calK(\bx,\bx') \delta c^{-2} (\bx') \dd \bx' \, ,  \label{def:IBS}
\end{align}
where $\widehat{\delta c} {}^{-2} (\bx) $ is the estimated velocity perturbation. We will use $\left< \cdot \right>$ to denote the statistical average with respect to the distribution of the random time delays or the stationary random source signals. The statistical average of the kernel $\calK$ (which is not the same kernel as the $\calK_\omega$ used earlier) is given by
\begin{align*} 
\langle \calK  (\bx,\bx') \rangle&=  \langle \calF^* \calF \rangle (\bx,\bx')    \nonumber  \\
&= \frac{1}{2\pi} \int \omega^4  |\hat{f} (\omega) |^2 \, \sum_{s=1}^{N_s}  \Big[ 
\overline{\hat{G}}\big(\omega, \bx,\by_s)  {\hat{G}}\big(\omega, \bx',\by_s \big)   \Big]     \sum_{r=1}^{N_r}   \Big[  \overline{\hat{G}}\big(\omega, \bx_r, \bx ) {\hat{G}}\big(\omega, \bx_r, \bx' \big)    \Big] \dd \omega   \nonumber \\
& = \frac{1}{2^5\pi^3} \int \omega^4  |\hat{f} (\omega) |^2 \sinc^2\Big(\frac{\omega}{c_0} |\bx-\bx'|\Big) \dd \omega   
\end{align*}
for noise blended sources, and by
\begin{equation*}
\langle \calK  (\bx,\bx') \rangle = \frac{T}{2^5\pi^3} \int \omega^4  \hat{F} (\omega) \sinc^2\Big(\frac{\omega}{c_0} |\bx-\bx'|\Big) \dd \omega
\end{equation*}
for the stationary random sources, see \cite{DFGS12}. 

We stress that the focus of this paper is on the high frequency analysis of this functional. We will analyze its performances in localizing singularities in the background velocity of propagation, in the high frequency regime $\eta\ll1$, where
\begin{equation*}
\frac{\omega_0}{c_0} = \frac{1}{\eta} \ .
\end{equation*}
Note that the wavelength is $2\pi\eta$.

\bigskip

\section{Expected contrast of the estimated perturbation}\label{sec:E}

Even if the estimated perturbation $\widehat{\delta c} {}^{-2}(\bx)$ provided by the imaging functional does not have the exact shape of the real perturbation, due to the different approximations used, it still shows a peak on the actual location of the velocity perturbation. We study in this section the expected (average in the statistical sense) contrast seen in the estimated velocity perturbation between the location of the real perturbation and points far from it.

The expected estimated perturbation for noise blended sources is given by 
\begin{align*}
\langle \widehat{\delta c} {}^{-2}(\bx)  \rangle  =  \frac{1}{2^5\pi^3} \int \calI_1(\bx,\omega) \omega^4 \big| \hat f (\omega) \big|^2 \dd \omega \ ,
\end{align*}
with
\begin{align*}
 \calI_1(\bx,\omega) = \int_{B_R} \sinc^2\Big(\frac{\omega}{c_0} |\bx-\bx'|\Big)  \delta c^{-2} (\bx') \dd \bx' \ .
\end{align*}
When the bandwidth of the source pulse is smaller than the central frequency, it is enough to study the behavior of the spatial integral $\calI_1(\bx)$ at the central frequency. For simplicity we drop the dependence on $\omega$. It turns out that $T \cdot \calI_1(\bx)$ is the quantity one has to study also in the case of stationary random sources. All the computations presented in the remaining part of this section are performed in the noise blended case; to obtain the different orders of amplitude provided by the imaging functional for stationary random sources it suffices to multiply the results by $T$.

To simplify the presentation of some computations, we assume that the support of the perturbation is centered on the center of the ball $B_R$ on the surface of which are located the sources and receivers, and choose this center as the origin of our system of coordinates.

For the purposes of this section it would not be necessary to distinguish between the scale of the wavelength of the signals ($\eta$) and that of the size of the perturbation ($\varepsilon$). However, this distinction turns out to be crucial in the analysis of fluctuations carried out in section \ref{sec:oscillation}.

\subsection{Point singularities}
Let us start by considering a point singularity, that is to say, a singularity whose support has a very small diameter. We model it with a perturbation of the velocity of propagation that is supported on a ball of radius $\varepsilon$:
\begin{equation} \label{eq:point}
\delta c^{-2} (\bx) = \alpha \, \indicator{B_\varepsilon} (\bx) \, .
\end{equation}
To simplify computations, we take $\alpha$ to be constant. Then, one observes that it only enters formulas as a multiplicative constant (squared in section \ref{sec:oscillation}): since it is of no relevance to our scopes, we set it equal to 1, also for the other types of perturbations.

Changing variables, we have:
\begin{align*}
\calI_1(\bx) &= \int_{B_R} \sinc^2 \big(|\bx-\bx'|/\eta \big) \delta c^{-2} (\bx')  \dd \bx' =  \int_{B_\varepsilon} \sinc^2 \big(|\bx-\bx'|/\eta \big) \dd \bx' \\
&=  \varepsilon^3 \int_{B_1} \sinc^2 \Big(\frac{\varepsilon}{\eta}| \bx/\varepsilon-  \bx'|\Big) \dd \bx' \ .
\end{align*}
At the center of the perturbation, $\bx=0$, we have
\begin{align*}
\calI_1(0) &=  \varepsilon^3 \int_{B_1} \sinc^2 \Big( \frac{\varepsilon}{\eta}|\bx'| \Big) \dd \bx' = 4\pi  \varepsilon^3 \int_0^1 \frac{\eta^2}{\varepsilon^2}  \sin^2 \Big( \frac{\varepsilon}{\eta} \rho \Big) \dd \rho =  2\pi \varepsilon \eta^2 \Big( 1- \sinc(2\varepsilon/\eta) \Big) .
\end{align*}
For $\varepsilon \simeq \eta$ we have that $\calI_1(0)=O(\varepsilon \eta^2)$. But if $\varepsilon \ll \eta$, the order becomes $O(\varepsilon^3)$.

For $|\bx| = O(1)$, what gives the order of amplitude of $\calI_1(\bx)$ is the decay of the $sinc(x)$ function, which goes approximately as $1/x$. We have that 
\begin{align*}
\calI_1(\bx) \simeq  \varepsilon^3 \int_{B_1}  \eta^2 |\bx-\varepsilon \bx'|^{-2} \, \sinc^2\big(|\bx|/\eta\big)  \dd \bx' \lesssim \frac{4\pi}{3} \varepsilon^3 \eta^2.
\end{align*}
Remark that the bound found is sharp, since there are $\bx$ for which $\calI_1$ is exactly of order $ \varepsilon^3 \eta^2$.

We see that the difference in amplitude between the centre of the perturbation and a point far from it is significant: it is of order $\eta^{-2}$ when $\varepsilon \ll \eta$, $\varepsilon^{-2}$ if $\varepsilon \simeq \eta$.

\subsection{Line singularities}
We consider now line--type singularities, which is to say an almost one--dimensional perturbation of the velocity of propagation. We model it by a perturbation supported on a cylinder of radius $\varepsilon$:
\begin{equation} \label{eq:line}
\delta c^{-2} (\bx) =  \indicator{C_\varepsilon} (\bx) \, , \qquad \qquad  C_\varepsilon = B_\varepsilon \times [-1,1] \subset \R^2\times \R \, .
\end{equation}
We have
\begin{align*}
\calI_1(\bx) &= \int_{B_R} \sinc^2 \big( |\bx-\bx'|/\eta \big) \delta c^{-2} (\bx')  \dd \bx' \\
&  =  \varepsilon^2 \int_{-1}^{1} \int_{B_1} \sinc^2 \Big(\frac{1}{\eta} \sqrt{(x-\varepsilon x')^2 + (y-\varepsilon y')^2 + (z-z')^2 }\Big) \dd \bx' \, .
\end{align*}
At $\bx=0$ this term reduces to
\begin{align*}
 \varepsilon^2 \int_{-1}^{1} \dd z' \int_{B_1} \sinc^2 \Big(\frac{1}{\eta} \sqrt{\varepsilon^2 (x'^2 + y'^2) + z'^2 }\Big) \dd x' \dd y'  \, ,
\end{align*}
which is rapidly oscillating in $z'$. In the integral in $z'$ we can therefore approximate $sin^2$ by its mean:
\begin{align*}
\calI_1(0) & \simeq 2\pi  \varepsilon^2 \int_{0}^1  \, \bigg( \frac{1}{2} \int_{-1}^{1} \eta^2 \Big(\varepsilon^2 r^2 + z'^2 \Big)^{-1}  \dd z'  \bigg)  r \dd r = 2\pi   \varepsilon^2 \int_0^1  \frac{\eta^2}{\varepsilon r}  \arctan\Big( \frac{1}{\varepsilon r}  \Big) \, r \dd r \\
&=  2\pi  \varepsilon \eta^2 \int_0^1 \arctan\Big( \frac{1}{\varepsilon r}  \Big) \dd r = O(\varepsilon \, \eta^2) \, .
\end{align*}
Far from the perturbation the integral is of order $\varepsilon^2\eta^2$. For example, for $x^2 + y^2 = C = O(1)$ we have
\begin{align*}
\calI_1(\bx) &\lesssim \frac{1}{2}   \varepsilon^2 \int_{-1}^{1} \int_{B_1} \eta^2 \Big( (x-\varepsilon x')^2 + (y-\varepsilon y')^2 + (z-z')^2 \Big)^{-1} \dd \bx' \\
&\simeq  \frac{\pi}{2} \varepsilon^2 \eta^2  \int_{-1}^{1}  \frac{1}{C+(z-z')^2} \dd z'   = O(\varepsilon^2 \, \eta^2) \, .
\end{align*}

Therefore, the difference in amplitude seen between the centre of the line perturbation and a point far from it is of order $\varepsilon^{-1}$.

\subsection{Plane singularities}\label{sec:speranza planes}
Let us consider now singularities that are approximately two--dimensional: we call them surface--type singularities and model them by a perturbation of the velocity which is supported on a disc of thickness $\varepsilon$:
\begin{equation}\label{eq:plane}
\delta c^{-2} (\bx) =  \indicator{D_\varepsilon} (\bx) \, , \qquad \qquad  D_\varepsilon = [-\varepsilon,\varepsilon] \times B_1 \subset \R\times \R^2 \, .
\end{equation}
With the notation introduced above we have
\begin{align*}
\calI_1(\bx) &  =  \varepsilon \int_{-1}^{1} \int_{B_1} \sinc^2 \Big(\frac{1}{\eta} \sqrt{(x-\varepsilon x')^2 + (y-y')^2 + (z-z')^2 }\Big) \dd \bx' \ .
\end{align*}
At $\bx=0$
\begin{align*}
\calI_1(0) &= 2\pi  \varepsilon \int_{-1}^{1} \dd x' \int_{0}^1 \sinc^2 \Big(\frac{1}{\eta} \sqrt{\varepsilon^2 x'^2 +r^2 }\Big) \ r \dd r \\
&\simeq  \pi  \varepsilon \int_{-1}^{1} \dd x' \int_{0}^1 \eta^2 \Big( \varepsilon^2 x'^2 +r^2 \Big)^{-1} \ r \dd r = \pi \varepsilon \eta^2    \int_{-1}^{1}  \frac{1}{2} \ln\Big( 1+ \frac{1}{\varepsilon^2 x'^2} \Big) \dd x' \\
&= O\big(\varepsilon \, \eta^2 |\ln(\varepsilon)| \big) 
\end{align*}
because the $sin^2$ in the first line is rapidly oscillating in $r$.

For $x = O(1)$ we have instead that 
\begin{align*}
\calI_1(\bx) &\lesssim  \varepsilon \int_{-1}^{1} \int_{B_1} \eta^2  \Big( (x-\varepsilon x')^2 + (y-y')^2 + (z-z')^2 \Big)^{-1} \dd y' \dd z' \dd x'   \\
&\simeq   2\pi   \varepsilon \eta^2 \int_{0}^{1} \frac{r}{x+r^2} \dd r  =  O(\varepsilon \, \eta^2)    \, .
\end{align*}
For this type of singularities, the difference in amplitude seen between the centre of the perturbation and a point far from it is very weak: it is only of order $| ln(\varepsilon)|$. This already hints to the fact that, with the imaging functional we consider, surface--type singularities are harder to locate than the other types.

\bigskip

\section{Fluctuations} \label{sec:oscillation}

We have obtained in the previous section the average contrast of the imaged perturbation. We want now to find confidence intervals for the typical contrast observed during the experiment. To do so, we must analyze the fluctuations in the result. They are given by the standard deviation of the estimated perturbation $ \widehat{\delta c}{}^{-2}$, which at a point $\bx$ is given by 
\begin{equation}\label{eq:Var}
\calS (\bx) = \Big \langle \, \big| \widehat{\delta c} {}^{-2} (\bx) - \langle   \widehat{\delta c} {}^{-2}  (\bx) \rangle \big|^2 \, \Big\rangle^\frac{1}{2} \ .
\end{equation}
We need to compute the standard deviation at the location of the perturbation and at points far from it, and compare them with the expected amplitude of the estimated perturbation. 
Using (\ref{def:IBS}) to write explicitly (\ref{eq:Var}) we get
\begin{align*}
\calS(\bx) &= \Big \langle \Big[ \int_{B_R} \big( \calK(\bx,\bx') - \langle\calK(\bx,\bx') \rangle \big) \delta c^{-2} (\bx') \dd \bx' \Big]^2 \Big \rangle^\frac{1}{2} \\
& = \bigg[ \iint_{B_R} \Big \langle \Big(  \calK(\bx,\bx') - \langle\calK(\bx,\bx') \Big)  \Big(  \calK(\bx,\bx'') - \langle\calK(\bx,\bx'') \Big) \Big \rangle \ \delta c^{-2} (\bx')  \delta c^{-2} (\bx'')  \dd \bx'  \dd \bx'' \bigg] ^\frac{1}{2} \\
&  = \bigg[ \iint_{B_R}  \mathrm{Cov} \Big(  \calK(\bx,\bx') ,  \calK(\bx,\bx'') \Big) \, \delta c^{-2} (\bx')  \delta c^{-2} (\bx'')  \dd \bx'  \dd \bx'' \bigg]^\frac{1}{2}  .
\end{align*}
In \cite{DFGS12} a formula was obtained for the variance of the kernel $\calK$. It is possible to carry out the same computations  for the covariance: for noise blended sources one obtains
\begin{align*}
\nonumber
& 2 \pi \, T_\tau \, {\rm Cov} \Big(  \calK(\bx,\bx') ,  \calK(\bx,\bx'') \Big) \simeq
 \int \dd \omega   |\hat{{f}}(\omega)|^4 \omega^8 \\
 & \times \bigg\{\sum_{r=1}^{N_r}   \overline{\hat{G}}\big(\omega, \bx_r, \bx ) \hat{G}\big(\omega,\bx_r, \bx' \big) \  \sum_{r=1}^{N_r}  \hat{G}\big(\omega, \bx_r, \bx ) \overline{\hat{G}}\big(\omega,\bx_r, \bx'' \big)    \nonumber \\
 \nonumber
  &  \qquad \times 
\Big\{ \sum_{s=1}^{N_s} |\hat{G}\big(\omega, \bx,\by_s) |^2  \  \sum_{s=1}^{N_s}  \hat{G}\big(\omega, \bx' ,\by_s)  \overline{\hat{G}}\big(\omega, \bx'' ,\by_s)   - \sum_{s=1}^{N_s} |\hat{G}\big(\omega, \bx,\by_s) |^2 \hat{G}\big(\omega, \bx' ,\by_s) \overline{\hat{G}}\big(\omega, \bx'' ,\by_s)    \Big\}
 \\
 \nonumber
& \quad +  \sum_{r=1}^{N_r} \overline{\hat{G}}\big(\omega, \bx_r, \bx ) \hat{G}\big(\omega,\bx_r, \bx' \big) \   \sum_{r=1}^{N_r} \overline{\hat{G}}\big(\omega, \bx_r, \bx ) \hat{G}\big(\omega,\bx_r, \bx'' \big)  \\
& \qquad \times  
\Big\{  \sum_{s=1}^{N_s} \overline{\hat{G}}\big(\omega, \bx,\by_s) {\hat{G}}\big(\omega, \bx',\by_s)    \  \sum_{s=1}^{N_s} \overline{\hat{G}}\big(\omega, \bx,\by_s) {\hat{G}}\big(\omega, \bx'',\by_s)       \nonumber \\
& \qquad \qquad  \qquad -\sum_{s=1}^{N_s} \overline{\hat{G}}\big(\omega, \bx,\by_s)^2 \hat{G}\big(\omega, \bx' ,\by_s)  \hat{G}\big(\omega, \bx' ,\by_s)  \Big\}  \ \bigg\}  \, .
\nonumber 
\end{align*}
Here, $T_\tau = \big( \int p_\tau^2 (t) \dd t \big)^{-1}$ and $p_\tau(t)$ is the probability density function of the random time delays. We will see that this quantity should be large so that the kernel $\calF^*\calF (\bx,\bx')$ is statistically stable. Using the method of Lagrange multipliers one can show for example that the maximal value of $T_\tau$ amongst all probability density functions $p_\tau$ compactly supported in $[-\tau_{max},\tau_{max}]$ is obtained for the uniform density over the interval and gives $T_\tau=2\tau_{max}$. Since $\tau_{max}$ must be at most of the order of the recording time $\tau_{max}\simeq T/2$, to obtain a large value of $T_\tau$ one should take $T$ large too (recall that in section \ref{sec:direct-inverse Pb} we had already assumed $T$ to be large). 

Observe that in the above equation for the covariance, in each of the two terms on the right hand side we are summing $N_s$ terms with a minus sign and $N_s^2$ with a plus sign, and they are all of the same order. The contribution of the terms with a minus sign is therefore small, and we can use the results of section \ref{sec:imaging func} to simplify the above equation into  \begin{align*}
 {\rm Cov} \Big(  \calK(\bx,\bx') ,  \calK(\bx,\bx'') \Big) \simeq \frac{1}{2 \pi \, T_\tau } \int & \dd \omega   |\hat{{f}}(\omega)|^4 \omega^8 \\
& \times  \bigg\{ \sinc\Big(\frac{\omega}{c_0} |\bx-\bx'|\Big)  \sinc\Big(\frac{\omega}{c_0} |\bx-\bx''|\Big)  \sinc\Big(\frac{\omega}{c_0} |\bx'-\bx''|\Big) \\
&   \hspace{1cm}  +  \sinc^2\Big(\frac{\omega}{c_0} |\bx-\bx'|\Big)  \sinc^2\Big(\frac{\omega}{c_0} |\bx-\bx''|\Big)  \bigg\} \, .
\end{align*}
Similar computations for the case of stationary random sources, where the terms with a minus sign do not even appear, leads to
\begin{align*}
{\rm Cov} \Big(  \calK(\bx,\bx') ,  \calK(\bx,\bx'') \Big) \simeq \frac{T}{2 \pi } \int & \dd \omega   |\hat{F}(\omega)|^2 \omega^8  \\
&  \times  \bigg\{ \sinc\Big(\frac{\omega}{c_0} |\bx-\bx'|\Big)  \sinc\Big(\frac{\omega}{c_0} |\bx-\bx''|\Big)  \sinc\Big(\frac{\omega}{c_0} |\bx'-\bx''|\Big) \\
&   \hspace{1cm}  +  \sinc^2\Big(\frac{\omega}{c_0} |\bx-\bx'|\Big)  \sinc^2\Big(\frac{\omega}{c_0} |\bx-\bx''|\Big)  \bigg\} \, .
\end{align*}

Let us continue the computations in the case of noise blended sources. Putting everything together, we find that the standard deviation of the imaging functional at a point $\bx$ is given by 
\begin{align*}
\calS(\bx)  = \bigg[ \frac{1}{2 \pi \, T_\tau } \int   |\hat{{f}}(\omega)|^4 \omega^8 \ \calI_2(\bx,\omega) \dd \omega  \bigg] ^\frac{1}{2}  
\end{align*}
with
\begin{align*}
\calI_2(\bx, \omega)&= \iint_{B_R}  \bigg[ \sinc\Big(\frac{\omega}{c_0} |\bx-\bx'|\Big)  \sinc\Big(\frac{\omega}{c_0} |\bx-\bx''|\Big)  \sinc\Big(\frac{\omega}{c_0} |\bx'-\bx''|\Big)   \\
&   \hspace{12mm}  +  \sinc^2\Big(\frac{\omega}{c_0} |\bx-\bx'|\Big)  \sinc^2\Big(\frac{\omega}{c_0} |\bx-\bx''|\Big)  \bigg]     \delta c^{-2} (\bx')  \delta c^{-2} (\bx'')  \dd \bx'  \dd \bx'' .  
\end{align*}

Again, we focus on the spatial integral $\calI_2$ (and to simplify notations we drop the dependence on $\omega$): the quantity we need to study is $\big[ \calI_2 (\bx) / T_\tau\big]^{1/2}$ for the noise blended sources and $\big[T \cdot \calI_2 (\bx)\big]^{1/2}$ for the stationary random sources. The computations presented below are performed in the noise blended sources setting. To obtain the corresponding standard deviation for stationary random sources it suffices to substitute the factor $1/T_\tau$ (or $1/\sqrt{T_\tau}$) by $T$ (or $\sqrt{T}$). For comparison with the average amplitude, recall that for stationary random sources the results obtained in the previous section have to be multiplied by a factor $T$.

We can rewrite the integral $\calI_2 (\bx)$ as the sum of the two integrals
\begin{align*}
\calJ_1 (\bx) &=  \iint_{B_R} \calH_1 (\bx,\bx',\bx'') \delta c^{-2} (\bx')\delta c^{-2} (\bx'') \  \dd \bx' \dd \bx'' \ , \\
\calJ_2 (\bx) &=  \iint_{B_R} \calH_2 (\bx,\bx',\bx'') \delta c^{-2} (\bx')\delta c^{-2} (\bx'') \  \dd \bx' \dd \bx'' \ ,
\end{align*}
 which is to say the double integral of the two kernels
\begin{align*}
\calH_1 (\bx,\bx',\bx'') &= \sinc\Big(\frac{\omega}{c_0} |\bx-\bx'|\Big)  \sinc\Big(\frac{\omega}{c_0} |\bx-\bx''|\Big)  \sinc\Big(\frac{\omega}{c_0} |\bx'-\bx''|\Big) \, ,\\
\calH_2 (\bx,\bx',\bx'') &= \sinc^2\Big(\frac{\omega}{c_0} |\bx-\bx'|\Big)  \sinc^2\Big(\frac{\omega}{c_0} |\bx-\bx''|\Big) 
\end{align*}
applied to the perturbation $\big( \delta c^{-2} (\bx')\delta c^{-2} (\bx'') \big)$, for every one of the three types of perturbation studied above. However, since
$$\calJ_2(\bx) = \calI_1^2(\bx)\, ,$$
we will only need to analyze $\calJ_1$.

For the following estimates it is important to keep separated the scale of the dimension of the perturbation $(\varepsilon)$ from the scale of the wavelength of the sources $(\eta)$. Their relative amplitude will be specified, but to help the reader to keep track of the different orders, we stress that we will always have $\varepsilon\le\eta$.

\subsection{Point singularities}
We return to the case of point perturbations introduced in the previous section and modeled by (\ref{eq:point}). In this case, even a very rude estimation is sufficient to obtain a bound which guarantees that this perturbation can be imaged. Since $|sinc| \le 1$, changing variables we get
\begin{align}
\big| \calJ_1(\bx) \big| &  \le \varepsilon^6 \iint_{B_1}\big|  \sinc\big( |\bx-\varepsilon \bx'|/\eta \big)  \sinc\big( |\bx- \varepsilon \bx''|/\eta \big) \big| \dd \bx' \dd \bx''       \nonumber  \\
&  \le K \varepsilon^6 \int_{B_1(\bx/\varepsilon)} \sinc^2 \big( |\bx'| \varepsilon/\eta \big) \dd \bx' \, ,  \label{eq:H1-points}
\end{align} 
where the constant $K=4/3 \pi$ comes from Cauchy-Schwarz inequality. Since the $sinc$ function is bounded, we have obtained a bound of order $\varepsilon^6$. This is only a rough upper bound, but it is not necessary to look for an improvement since it is already of the same order of the integral of the second kernel, for which we have
\begin{align*}
\calI_2(0) = \calI_1^2 (0) \simeq \varepsilon^6 \, .
\end{align*}
Using the decay of the $sinc$ function, we can find near the perturbation ($|\bx| \simeq\varepsilon$) a bound for $\calI_2$ of the same order. Oscillations are therefore of order $\varepsilon^3/\sqrt{T_\tau}$.

Recall that the average value observed on the peak is of order $\varepsilon^3$ for $\varepsilon \ll \eta$, so that the typical value observed remains of the same order due to the large factor $T_\tau$. The same results holds true also for $\varepsilon \sim \eta$, namely the typical value observed is of the same order of the average value. \\

Far from the perturbation the integrals of the two kernels decrease. Using the bound (\ref{eq:H1-points}), the integral of the first kernel can be bounded for $|\bx|= O(1)$ by
\begin{equation*}
\big| \calJ_1(\bx)  \big| \le K \varepsilon^6 \int_{B_1} \eta^2 \big| \bx - \varepsilon \bx' \big|^{-2} \dd \bx' =  K \, \varepsilon^6 \eta^2 \, .
\end{equation*}
With some more work, one can show that this bound is sharp, at least for $\varepsilon \ll \eta$. This can be done using the Fourier representation of the $sinc$ function written in spherical coordinates $\bu = \bu(r,\theta,\phi) \in \R^3$:
\begin{align}
\sinc( \lambda|\bx| ) &=  \frac{1}{2} \int_{-1}^{1} e^{-i \lambda\zeta |\bx|} \dd \zeta = \frac{1}{2} \int_{0}^{\pi} e^{-i \lambda |\bx| \cos(\theta) } \sin (\theta)\dd \theta    \nonumber \\
&= \frac{1}{4\pi} \int_{0}^{2\pi} \int_0^{\pi} e^{-i \lambda \bx \cdot \bu(r,\theta,\phi)} \sin (\theta)\dd \theta \dd \phi = \frac{1}{4\pi} \int_{S^2} e^{i \lambda \bx \cdot \bu} \dd \bu \ ,  \label{eq:sinc->exp}
\end{align}
where $S^2=\{\bx\in\R^3 : |\bx|=1 \} $ is the unitary sphere in $\R^3$. We can write
\begin{align*}
\calJ_1 (\bx) &= (4\pi)^{-3} \iint_{B_\varepsilon} \iiint_{S^2} e^{i\frac{1}{\eta} (\bx-\bx')\cdot \bu} e^{i\frac{1}{\eta} (\bx-\bx'')\cdot \bv} e^{i\frac{1}{\eta} (\bx'-\bx'')\cdot \bw} \dd \bu \dd \bv \dd \bw \ \dd \bx' \dd \bx''  \nonumber \\
&\simeq  \frac{\varepsilon^6}{3^2 4 \pi}  \iiint_{S^2} e^{i\frac{1}{\eta} \bx \cdot (\bu + \bv)}  \dd \bu \dd \bv \dd \bw \ , 
\end{align*}
where we have used the assumption $\varepsilon \ll \eta$. Simplifying this equation and using again (\ref{eq:sinc->exp}), we get for $|\bx|=O(1)$
\begin{align*}
\calJ_1 (\bx) &\simeq \frac{\varepsilon^6}{3^2}   \Big[ \int_{S^2} e^{i\frac{1}{\eta} \bx \cdot \bu}  \dd \bu \Big]^2 = \frac{\varepsilon^6}{3^2}   \Big[ 4 \pi \, \sinc\big( |\bx| /\eta \big)  \Big]^2 = O(\varepsilon^6 \eta^2) \ .  
\end{align*}
As for the integral of the second kernel, the bound we have is of order $\varepsilon^{6}\eta^4$. For $\varepsilon\ll \eta$ we have therefore
\begin{align*}
\calI_2(\bx) = O(\varepsilon^6\eta^2) \, .
\end{align*}
In the general case $\varepsilon\le\eta$, the above equation becomes an upper bound.

Recall that the statistical average of the imaging operator far from the perturbation is of order $\varepsilon^3\eta^2$. Assuming that $T$ is large, but still $1\ll T_\tau \le 1/\eta^2$, the above result implies that the typical value observed is at most of order $\varepsilon^3\eta/\sqrt{T_\tau}$ (it is exactly of this order for $\varepsilon \ll \eta$). Therefore, the typical contrast is still at least of order $ \sqrt{T_\tau}/ \eta$, allowing for a precise location of the perturbation (both when $\varepsilon\ll\eta$ and $\varepsilon\sim \eta$).

\subsection{Line singularities}
Consider the case of line singularities, modeled by (\ref{eq:line}). Using rude estimations similar to the ones presented above, we could only bound the integral of the (absolute value of the) first kernel near the origin with something of order $\varepsilon^4\eta^2\ln^2(\varepsilon)$. This means that, in order to be sure to able to image the perturbation, we would need to have $\varepsilon |\ln(\varepsilon)| \ll \eta$. But we can do better.\\

Assuming simply $\varepsilon\ll\eta$, we can approximate  
\begin{equation}\label{eq:approx sinc}
\sinc \Big(\frac{1}{\eta} \sqrt{\varepsilon^2( x'^2 + y'^2) + z'^2 }\Big) \simeq \sinc \big( |z'|/\eta \big) = \sinc \big( z'/\eta \big) \ .
\end{equation}
It is then possible to use the Fourier representation of the $sinc$ function
\begin{equation*}
\sinc(\lambda z) = \frac{1}{2\lambda} \int_{-\lambda}^{\lambda} e^{-i \zeta z} \dd \zeta =\frac{1}{2} \int_{-1}^{1} e^{-i \lambda\zeta z} \dd \zeta
\end{equation*}
to obtain the amplitude of oscillations. At $\bx=0$ we can rewrite the integral of the first kernel as an integral over $C_1 = B_1 \times [-1,1] \subset \R^2\times \R$, use (\ref{eq:approx sinc}) and integrate in $x', y', x'', y''$:
\begin{align*}
\calJ_1(0) &=  \iint_{C_\varepsilon} \sinc \big( |\bx'| /\eta \big)  \sinc\big( |\bx''| /\eta \big)  \sinc\big( |\bx'-\bx''| /\eta \big) \dd \bx' \dd \bx''\\
&=  \varepsilon^4 \iint_{C_1} \sinc \Big(\frac{1}{\eta}  \sqrt{\varepsilon^2( x'^2 + y'^2) + z'^2 }\Big)   \sinc \Big(\frac{1}{\eta}  \sqrt{\varepsilon^2( x''^2 + y''^2) + z''^2 }\Big) \times\\
& \hspace{2cm} \times \sinc \Big(\frac{1}{\eta}  \sqrt{\varepsilon^2\big( (x'-x'')^2 + (y'-y'')^2\big) + (z'-z'')^2 }\Big)    \dd \bx' \dd \bx'' \\
&\simeq \pi^2  \varepsilon^4 \iint_{z-1}^{z+1} \sinc \big( z'/\eta \big)\sinc \big( z''/\eta \big)  \sinc \big( (z'-z'')/\eta \big) \dd z' \dd z'' \ ,
\end{align*}
and using the Fourier representation introduced above we get 
\begin{align*}
\calJ_1(0) &= \frac{\pi^2}{2^3}  \varepsilon^4 \iint_{z-1}^{z+1} \iiint_{-1}^1 e^{-i  \zeta' z'/\eta} e^{-i  \zeta'' z''/\eta} e^{-i  \zeta (z'-z'')/\eta} \dd \zeta \dd \zeta' \dd \zeta'' \dd z' \dd z''\\
&=  \frac{\pi^2}{2}\,  \varepsilon^4  \iiint_{-1}^1 \sinc \Big( \frac{1}{\eta}(\zeta'+\zeta) \Big) \sinc \Big( \frac{1}{\eta} (\zeta''-\zeta) \Big)  \dd \zeta \dd \zeta' \dd \zeta'' \\
&=   \frac{\pi^2 }{2}\,  \varepsilon^4 \eta^2  \int_{-1}^1 \int_{(\zeta-1)/\eta}^{(\zeta+1)/\eta} \sinc(u_1)  \dd u_1 \int_{(-\zeta-1)/\eta}^{(-\zeta+1)/\eta}  \sinc(u_2)   \dd u_2   \ \dd \zeta \ .
\end{align*}
Since the function $s \mapsto \int_0^s \sinc(u) \dd u$ is uniformly bounded in $s$, we get that $$\calJ_1 (0) = O(\varepsilon^4 \eta^2)\ .$$

For the second kernel, we have
\begin{equation*}
\calJ_2(0) = \calI_1^2 (0) =  O(\varepsilon^2\eta^4)\, .
\end{equation*}


Remark that for $|\bx|\simeq \varepsilon$ we can still bound fluctuations in the same way, because (\ref{eq:approx sinc}) still holds, and the integral
\begin{align*}
\int \Big(  (x-\varepsilon x')^2 + (y-\varepsilon y')^2 + (z-z')^2 \Big)^{-1/2} \dd \bx'
\end{align*}
is maximal when $\bx=0$.

We see that fluctuations near $\bx=0$ are of order $\varepsilon \eta^2/\sqrt{T_\tau}$. Since the average value observed at $\bx=0$ for the imaging functional is of order $\varepsilon \eta^2$, it is thanks to the large factor $T_\tau \gg1$ that we get the statistical stability of the operator. This means that the typical value observed on the perturbation remains of order $\varepsilon\eta^2$.\\

When $\bx$ is far from the perturbation, oscillations are even smaller. Indeed, we can proceed as in (\ref{eq:H1-points}) to find a  bound for the integral of the absolute value of the first kernel. We get
\begin{align*}
\calJ_1(\bx) &\le   \varepsilon^4 \bigg[ \int_{-1}^1 \int_{B_1}\big| \sinc \Big(\frac{1}{\eta} \sqrt{(x-\varepsilon x')^2 + (y-\varepsilon y')^2 + (z-z')^2 }\Big)   \big| \dd x' \dd y' \ \dd z' \bigg]^2 \\
& \simeq   \varepsilon^4 \eta^2 \bigg[ \int_{-1}^1 \int_{B_1} \big( (x-\varepsilon x')^2 + (y-\varepsilon y')^2  + (z-z)'^2 \big)^{-1/2} \, \dd \bx' \bigg]^2  .
\end{align*}
Since the integrand is bounded, the bound we get is of order $\varepsilon^4\eta^2$. The second kernel is of a higher order, $\calJ_2(\bx) \lesssim O(\varepsilon^4\eta^4)$. 
This bound implies that the typical value observed far from the perturbation is of order at most $\varepsilon^2\eta/\sqrt{T_\tau}$, so that the contrast is at least of order $\eta\sqrt{T_\tau}/\varepsilon$.

\subsection{Plane singularities}
We turn now to analyze fluctuations in the case of surface--type perturbations, modeled by (\ref{eq:plane}). The difficult part is again to obtain good estimates on the integral of the first kernel, for which we use the Fourier representation of the $sinc$ function obtained in (\ref{eq:sinc->exp}). At $\bx=0$ we have
\begin{align*}
\calJ_1(0) &= \frac{1}{(4\pi)^3}\iint_{D_\varepsilon} \iiint_{S^2} e^{i\frac{1}{\eta} [ \bx' \cdot (\bu+\bw)  + \bx'' \cdot(\bv-\bw) ] } \dd \bu \dd \bv \dd \bw \ \dd \bx' \dd \bx''\\
&= \frac{1}{(4\pi)^3} \int_{S^2} \bigg[\int_{S^2} \int_{D_\varepsilon}  e^{i\frac{1}{\eta} \bx' \cdot (\bu+\bw) } \dd \bx' \ \dd \bu \bigg]^2 \dd \bw \\
&\simeq \frac{4\varepsilon^2}{(4\pi)^3} \int_{S^2} \bigg[ \int_{S^2} \int_{B_1} e^{i\frac{1}{\eta}  \bx'_\perp \cdot (\bu+\bw)  } \dd \bx'_\perp  \dd \bu \bigg]^2 \dd \bw \ ,
\end{align*}
where the approximate equality holds for $\varepsilon\ll \eta$ and $\perp$ denotes the projection on the last two coordinates: $\bx_\perp = (y,z) \in \R^2$. The integral in $\dd \bx'_\perp$ is computed on $(D_\varepsilon)_\perp = B_1\in\R^2$. We can rewrite also the integrals on $S^2$ as (twice the) integrals on the projection $B_1\in\R^2$, compute the integral in $\dd \bx'_\perp$ and change variables:
\begin{align*}
\calJ_1(0) & \simeq \frac{32\varepsilon^2}{(4\pi)^3} \int_{B_1} \bigg[ \iint_{B_1} e^{i\frac{1}{\eta}  \bx'_\perp \cdot (\bu+\bw)_\perp }  \sqrt{1-|\bu_\perp|^2} \dd \bx'_\perp   \dd \bu_\perp \bigg]^2   \sqrt{1-|\bw_\perp|^2}  \dd \bw_\perp \\  
&= \frac{2\varepsilon^2}{\pi}  \int_{B_1} \bigg[ \int_{B_1} \frac{J_1\big( |\bu_\perp+\bw_\perp|\eta \big)}{|\bu_\perp+\bw_\perp|/\eta}  \sqrt{1-|\bu_\perp|^2}  \dd \bu_\perp \bigg]^2 \sqrt{1-|\bw_\perp|^2}  \dd \bw_\perp   \nonumber  \\
&= \frac{2\varepsilon^2\eta^6}{\pi } \int_{B_{1/\eta}}  \bigg[  \int_{B_{1/\eta}(\bw_\perp)}  \frac{J_1( |\bu_\perp|)}{ |\bu_\perp|} \sqrt{1-|\bu_\perp-\bw_\perp|^2 \,\eta^2} \, \dd \bu_\perp \bigg]^2     \nonumber  \\
& \hspace{7,5cm} \times \sqrt{1-|\bw_\perp|^2 \,\eta^2} \, \dd \bw_\perp \ .   \nonumber  
\end{align*}
Here, $J_1$ is the Bessel function of the first kind. Let us focus on the integral inside the square brackets. Observe that the origin of our system of coordinates is always inside $B_{1/\eta}(\bw_\perp)$. Changing to polar coordinates we have
\begin{align*}
\calY &= \int_{B_{1/\eta}(\bw_\perp)} \frac{J_1( |\bu_\perp|)}{ |\bu_\perp|} \sqrt{1-|\bu_\perp-\bw_\perp|^2 \, \eta^2} \, \dd \bu_\perp  = \int_0^{2\pi}  \hspace{-2mm} \int_0^{\rho_\bw(\theta)} \hspace{-5mm}J_1(r) \phi_\bw(\theta,r) \dd r \dd \theta\ ,
\end{align*}
where we denote by $\phi_\bw$ the square root term (written in polar coordinates), and the function $\rho_\bw(\theta)$ takes values in $[1/\eta-|\bw_\perp|,1/\eta+|\bw_\perp|]\subset[0,2/\eta]$. We claim that the integral term $\calY$ is bounded. This can be proved integrating by parts in $r$.
Remark that the square root term is concave (as a function of $\bu_\perp$), take its maximum over the domain of integration $B_{1/\eta}(\bw_\perp)$ at the center of the ball and is zero at the boundary. Therefore, for every fixed $(\theta,\bw)$, the function $\phi_\bw(\theta,\cdot)$ is still concave and bounded by $1$. Then, one easily obtains that the integral of the absolute value of its derivative in $r$ is bounded by 2. Also, the antiderivative of $J_1(r)$ is the Bessel function of order zero $-J_0(r)$, which is bounded (the maximum of its absolute value is taken at $r=0$, and $J_0(0)=1$). Putting everything together, we get
\begin{align*}
\calY =& \int_0^{2\pi}  -J_0(r) \phi_\bw(\theta,r) \Big|_{r=0}^{r=\rho_\bw(\theta)}  + \int_0^{\rho_\bw(\theta)} J_0(r) \, \partial_r \phi_\bw(\theta,r) \dd r \ \dd \theta   \label{eq:calY piani}  \\
\le & \int_0^{2\pi} J_0(0) +  J_0(0) \int_0^{\rho_\bw(\theta)} \big| \partial_r \phi_\bw(\theta,r) \big| \dd r \ \dd \theta \  \le \, 6 \pi  \ .  \nonumber
\end{align*}
This proves the claim.

Using again the boundedness of the square root, we can bound the integral in $\dd \bw_\perp$ by $\pi/\eta^2$. We have therefore obtained a bound for $\calJ_1(0)$ of order $O(\varepsilon^2\eta^4)$. This is only an upper bound, but there is no need to look for an improvement, since it is already of a smaller order than the integral of the second kernel, for which we have
\begin{equation*}
\calJ_2(0) = \calI_1^2(0) = O(\varepsilon^2\eta^4 \ln^2(\varepsilon)) \, .
\end{equation*}
Therefore,
\begin{equation*}
\calI_2(0) \simeq \varepsilon^2 \eta^4 \ln^2(\varepsilon) \ .
\end{equation*}\\

Far from the perturbation, oscillations are even smaller. Denote $\bx= (x,\bx_\perp) \in\R\times\R^2$ and $\bu=(u_1,\bu_\perp) \in\R\times\R^2$; the same notation is used for $\bv$ and $\bw$. Let us look at the integral of the first kernel; following the computations presented above we have
\begin{align*}
\calJ_1(\bx) & = \frac{1}{(4\pi)^3}\hspace{-1mm} \iint_{D_\varepsilon} \iiint_{S^2} \hspace{-1mm} e^{i\frac{1}{\eta} [ (\bx-\bx') \cdot (\bu+\bw)  + (\bx-\bx'') \cdot(\bv-\bw) ] } \dd \bu \dd \bv \dd \bw \ \dd \bx' \dd \bx''    \\
&\simeq \frac{4\varepsilon^2}{(4\pi)^3} \iiint_{S^2} \iint_{B_1}  \hspace{-1mm}  e^{i\frac{1}{\eta} [ (\bx-\bx')_\perp^{} \cdot (\bu+\bw)_\perp^{}  + (\bx-\bx'')_\perp^{} \cdot(\bv-\bw)_\perp^{} ] } \dd \bx'_\perp \dd \bx''_\perp     \nonumber  \\
 &  \hspace{22mm} \times  e^{i\frac{1}{\eta}  [ x (u+w)_1^{}  + x (v-w)_1^{} ] } \dd \bu \dd \bv \dd \bw     \nonumber \\
&= \frac{32\varepsilon^2}{(4\pi)^3} \iiint_{B_1} \iint_{B_1}   \hspace{-1mm}  e^{-i\frac{1}{\eta} [ \bx'_\perp \cdot (\bu+\bw)_\perp^{}  + \bx''_\perp \cdot(\bv-\bw)_\perp^{} ] } \dd \bx'_\perp \dd \bx''_\perp \  \times e^{i\frac{1}{\eta} [ \bx_\perp^{} \cdot (\bu+\bv)_\perp^{}   ] }       \nonumber  \\
& \hspace{22mm}  \times e^{i\frac{1}{\eta}  x \big[ \sqrt{1- |\bu_\perp^{}|^2} + \sqrt{1 - |\bv_\perp^{}|^2} \big] }    \nonumber  \\
& \hspace{22mm}  \times \sqrt{1-|\bu_\perp|^2} \sqrt{1-|\bv_\perp|^2}\sqrt{1-|\bw_\perp|^2} \, \dd \bu_\perp \dd \bv_\perp \dd \bw_\perp    \nonumber  
\end{align*}
so that
 \begin{align*}
\calJ_1(\bx)&  \simeq \frac{2\varepsilon^2\eta^6}{\pi} \int_{B_{1/\eta}}  \hspace{-1mm}   \sqrt{1-|\bw_\perp|^2 \,\eta^2}  \nonumber   \\
&  \hspace{15mm}  \times  \hspace{-1mm}   \int_{B_{1/\eta}(\bw_\perp)} \hspace{-1cm} e^{i \bx_\perp^{} \cdot \, \bu_\perp^{} } e^{i x \sqrt{\frac{1}{\eta^2}  - \bu_\perp^2}  } \sqrt{1-|\bu_\perp-\bw_\perp|^2 \,\eta^2} \, \frac{J_1( |\bu_\perp|)}{ |\bu_\perp|} \dd \bu_\perp     \nonumber  \\
&  \hspace{15mm}  \times  \hspace{-1mm}  \int_{B_{1/\eta}(-\bw_\perp)} \hspace{-12mm} e^{i \bx_\perp^{} \cdot \, \bv_\perp^{} }  e^{i  x \sqrt{\frac{1}{\eta^2} - \bv_\perp^2}  } \sqrt{1-|\bv_\perp+\bw_\perp|^2 \,\eta^2} \, \frac{J_1( |\bv_\perp|)}{ |\bv_\perp|} \dd \bv_\perp \  \dd \bw_\perp \  .    \nonumber  
\end{align*}
For $|\bx|\gg1$, the last two integrals above are now much smaller than the corresponding ones for $\bx=0$. This is due to the fact that for $|\bx|$ large, at least one of the two exponential terms, which have mean zero, is rapidly oscillating with respect to $J_1$. We therefore have that $\calJ_1(\bx)$ is at most of order $\varepsilon^2\eta^4$. Far from the perturbation, the (sharp) bound we have on the integral of the second kernel is of the same order: $\calJ_2(\bx) \lesssim O(\varepsilon^2\eta^4)$.
 We have obtained that
\begin{equation*}
\calI_2(\bx) \simeq \varepsilon^2 \eta^4 \, .
\end{equation*}
Thanks to the large factor $T_\tau$, fluctuations are therefore smaller than the average value given by the imaging functional, both on the perturbation and far from it. The typical contrast remains therefore of the same order as the average contrast, namely of order $|\ln(\varepsilon)|$.

\bigskip

\section{Conclusions and comments}\label{sec:end}

We have analyzed the imaging functional 
 given by (\ref{eq:def F*})
 in the high frequency regime ($\eta\ll1$) with respect to small perturbations ($\varepsilon\ll1$) of the background velocity of propagation.
Using a suitable disposition of the sources and receivers, we have been able to obtain quantitative estimates on the (average, with respect to the realization of the random time delays or the stationary random source signals) sensitivity of the imaging functional. The image presents a peak on the location of the perturbation, and the contrast is of order $\eta^{-2}$ for point perturbations, of order $\varepsilon^{-1}$ for line perturbations, and only of order $|\ln(\varepsilon)|$ for surface perturbations.

The most interesting result obtained in this paper concerns the quantitative analysis of the statistical stability of this functional, providing the typical contrast seen for the three perturbations considered. The question of stability of the imaging functional has been addressed in \cite{DFGS12}: no quantitative analysis was carried out there, but it was shown that a condition for the statistical stability is that the quantities $T$ (for stationary random sources) and $T_\tau$ (for noise blended sources) must be large. For random time delays uniformly distributed on the interval $[-\tau_{max},\tau_{max}]$, $T_\tau$ large means that $\tau_{max}$ must be large, which in turn implies that the recording time $T\simeq\tau_{max}$ must be large. 

An important fact is that the typical contrast found only depends on the type of perturbation one is trying to image, and not on the method used. All results are described below for noise blended sources, but the corresponding contrast for stationary random sources are obtained simply substituting $T_\tau$ with $T$.
 
We have shown that for point perturbations, both when $\varepsilon\ll\eta$ and $\varepsilon\sim\eta$, fluctuations due to the stochastic nature of the method are small, and the typical contrast is at least of order $\sqrt{T_\tau}/\eta$ ($\sqrt{T}/\eta $ for stationary random sources): point perturbations are easy to find.

For line perturbation the situation is different. We can image with a satisfactory accuracy and reasonable recording time $T\gg1$ only very thin line perturbations, $\varepsilon\ll \eta$. The typical contrast in this case is at least of order $ \sqrt{T_\tau} \eta/\varepsilon$.

For plane perturbations the average contrast is quite poor, only of order $|\ln(\varepsilon)|$. However, for very thin perturbations  , $\varepsilon\ll \eta$, also the typical contrast is of the same order.

These results are summarized in the following tables, where we list for the three type of perturbations considered the order of the average value given by the imaging functional and of the standard deviation at the center of the perturbation and far from it. \\

\begin{table}[ht]
\centering
\begin{tabular} { | c | c c | c c |  }      
\hline \hline & & & &   \\ 
 &  $\big\langle \widehat{ \delta c} {}^{-2} (\bx)\big\rangle$ &   & $\calS(\bx)$ &  \\
 & $\bx = 0$  & $|\bx| \gg1$ & $\bx = 0$  & $|\bx| \gg1$   \\
 & & & &   \\
Points & $\simeq\varepsilon^3$ & $\lesssim \varepsilon^3\eta^2$ & $\simeq \varepsilon^3/\sqrt{T_\tau}$ &  $\lesssim \varepsilon^3\eta/\sqrt{T_\tau}$    \\
Lines & $\simeq\varepsilon\eta^2$  & $\lesssim \varepsilon^2\eta^2$ & $\simeq \varepsilon\eta^2 / \sqrt{T_\tau} $ & $\lesssim \varepsilon^2\eta/\sqrt{T_\tau}$   \\
Planes &$\simeq  \varepsilon\eta^2 |\ln(\varepsilon)|$ & $\lesssim  \varepsilon\eta^2$ & $\simeq \varepsilon\eta^2 |\ln(\varepsilon)| / \sqrt{T_\tau}$ & $\lesssim \varepsilon\eta^2/\sqrt{T_\tau}$  \\& & & &    \\
\hline 
\end{tabular}\vspace{3mm}
\caption { Noise blended sources: mean $\langle \widehat{ \delta c} {}^{-2} \rangle$ and standard deviation $(\calS)$ of the estimated velocity perturbation at the center of the perturbation $(\bx=0)$ and far from it $(|\bx|\gg1)$, in the regime $\varepsilon\ll \eta \ll 1$. The cases of point, line and disc singularities are displayed.}
\end{table}

\begin{table}[ht]
\centering
\begin{tabular} { | c | c c | c c |  }  
\hline \hline & & & &   \\ 
 &  $\big\langle \widehat{ \delta c} {}^{-2} (\bx)\big\rangle$ &   & $\calS(\bx)$ &  \\
 & $\bx = 0$  & $|\bx| \gg1$ & $\bx = 0$  & $|\bx| \gg1$   \\
 & & & &   \\
Points & $\simeq T \varepsilon^3$ & $\lesssim T \varepsilon^3\eta^2$ & $\simeq \sqrt{T} \varepsilon^3$ &  $\lesssim \sqrt{T} \varepsilon^3\eta$    \\
Lines & $\simeq T \varepsilon\eta^2$  & $\lesssim T \varepsilon^2\eta^2$ & $\simeq \sqrt{T} \varepsilon\eta^2  $ & $\lesssim \sqrt{T} \varepsilon^2\eta $   \\
Planes &$\simeq  T \varepsilon\eta^2 |\ln(\varepsilon)|$ & $\lesssim T  \varepsilon\eta^2$ & $\simeq \sqrt{T} \varepsilon\eta^2 |\ln(\varepsilon) | $ & $\lesssim \sqrt{T} \varepsilon\eta^2$ \\& & & &    \\
\hline
\end{tabular}\vspace{3mm}
\caption { Stationary random sources: mean $\langle \widehat{ \delta c} {}^{-2} (\bx)\rangle$ and standard deviation $(\calS)$ of the estimated velocity perturbation at the center of the perturbation $(\bx=0)$ and far from it $(|\bx|\gg1)$, in the regime $\varepsilon\ll \eta \ll 1$. The cases of point, line and disc singularities are displayed. }
\end{table}
 \vspace{5cm}

\section*{Acknowledgments} The author gratefully acknowledge the helpful and substantive critical comments by Prof. Josselin Garnier which greatly helped the development of this work. \\


\end{document}